\documentclass[onecolumn]{article}          
\usepackage{graphicx}
\usepackage{amssymb}
\usepackage{xypic}
\usepackage[english]{babel}
\usepackage{hyperref}

\newtheorem{theorem}{Theorem}[section]
\newtheorem{proposition}{Proposition}[section]
\newtheorem{lemma}[theorem]{Lemma}
\newtheorem{remark}{\it Remark\/}

\newcommand{\fiz}{{\textstyle{\frac{1}{2}}}}
\newcommand{\bA}{{\bf A}}
\newcommand{\bb}{{\bf b}}

\newcommand{\bP}{{\bf P}}

\newcommand{\bz}{{\bf z}}

\newcommand{\mL}{{\mathcal L} }
\newcommand{\mI}{{\mathcal I} }
\newcommand{\mJ}{{\mathcal J} }
\newcommand{\mK}{{\mathcal K} }

\begin{document}

\title{Towards exact symplectic integrators from Liouvillian forms
\thanks{The author was supported by a grant from the \emph{Fondation du Coll\`ege de France}
under the research convention PU14150472.}
}

\author{Hugo Jim\'enez-P\'erez
}

\maketitle

\begin{abstract}
In this article we introduce a low order implicit symplectic 
integrator designed to follow the Hamiltonian flow as close as 
possible. This integrator is obtained by the method of 
Liouvillian forms
and does not require particular hypotheses on the Hamiltonian.

The numerical scheme introduced in this paper is a 
modification of the symplectic mid-point rule,
it is symmetric and it is obtained by an isotopy of the 
deformation of the exact Hamiltonian flow to the straight line 
passing by two consecutive points of the discretized flow. This 
isotopy generates an alternative vector field on the flow lines
transversal to the Hamiltonian vector field. We consider only the 
line arising from the mid-point to construct the symplectic 
integrator.

\end{abstract}

\section{Introduction}
\label{}

In \cite{Jim15a} the author has introduced an alternative method for 
constructing implicit symplectic integrators using \emph{special symplectic manifolds}
\cite{Tul76,Tul77} and \emph{Liouvillian forms} \cite{LM87,Lib00}.
Such a method extends in a natural way 
the method of \emph{generating functions}, first introduced by Hamilton when working 
with optical paths \cite{Ham33} and then developed by Jacobi in \cite{Jac84}.
In a different context Poincar\'e developed the theory of integral invariants in 
his celebre \emph{Les m\'ethodes nouvelles de la m\'ecanique c\'eleste} \cite{Poi99}
where he used generating functions for studying bifurcating orbits arising from 
prescribed periodic orbits.
Generating functions were studied in symplectic geometry 
by many authors such as Viterbo \cite{Vit92}, Chaperon \cite{Cha95}
Maslov \cite{Mas65},  H\"ormander \cite{Hor71}, Weinstein \cite{Wei72} among many others.
From the numerical point of view Feng Kang and his coworkers 
\cite{Kang85a,KG88,GD95} have studied 
systematically the construction of symplectic integrators using generating
functions.  However, their point of view follows the Siegel's approach 
\cite{Sie64} which is based on the matrix algebra of the symplectic group.
A compilation of their work is contained in \cite{KQ10}. 

The relation between Liouvillian forms and generating functions is as follows. Using the 
Hodge decomposition of differential forms a Liouvillian form is decomposed in an exact, 
a harmonic and a co-exact 1-forms; this decomposition is unique \cite{Mor01}.
The exact part is related with the differential of a generating function, and 
they coincide on the Lagrangian surface defined by the generating function. 
The main difference between both methods is that a Liouvillian form is 
defined on open subsets of the symplectic manifold and it contains more information 
about its geometry than the generating function, among other advantages. 

For a problem with $n$ degrees of freedom a $n(n+1)$ dimensional 
continuous family of implicit symplectic 
integrators can be constructed under this method.
This was already noted by Kang and his coworkers \cite{KQ10}, however 
no geometric explanation concerning this family was given by them. 
In contrast, they interpret the Euler
symplectic methods as a first order approximation and the mid-point rule as
a second order approximation for the elements of this family of implicit 
symplectic maps. 
The method of Liouvillian forms gives a precise meaning to this family, a geometric 
explanation and a way to find an adapted symplectic integrator for a given 
(classical and natural) Hamiltonian problem.
The generating functions of type \emph{II,III} in \cite{Arn89}, 
(alternatively of type $V$ in \cite{MS17}), and the mid-point rules 
are just 3 different elements in the family. 
However, the generating functions of type \emph{I, IV} (alternatively of type $S$)
do not belog to this family. Moreover, the differential of the so called 
Poincar\'e's generating function \cite{Poi99}, which has been associated to the mid-point
rule, is a generating function for solving a different variational problem \cite{Jim15f}.

In the  \emph{method of Liouvillian forms}, the resolution of the 
Hamilton-Jacobi equation is not necessary and the algorithm is obtained from 
a suitable projection of the tangent space of a 2n-dimensional submanifold of 
the product of two symplectic manifolds, which is a Lagrangian submanifold with 
respect to the usual symplectic form. This submanifold is determined in a unique 
way by a triplet of Liouvillian forms.
The first numerical tests were shown in 
\cite{Jim15b}, where some Liouvillian forms were constructed in a random way.
At this point, the method has been completely formalized using differential
geometry. It lets us controlling the numerical solution since for every 
Liouvillian form we have, generically, a different integrator. 
In particular, we can control the oscillations of the numerical 
solution around the fixed value of the energy
and our interest becomes the search for the integrator which produces the minimal error. 
Liouvillian forms which are good candidates for 
integrators minimizing these oscillations, are close to those 
which produce the symplectic mid-point rule \cite{Jim15b}. 
Following the numerical evidence which predicts that the variation depends on the Hamiltonian,
we proved a series of results which explains this fact
\cite{Jim15b,Jim16a,JVR17,Jim18,Jim19a}.

\begin{figure}[h!]
 \centering
 \includegraphics[scale=0.20]{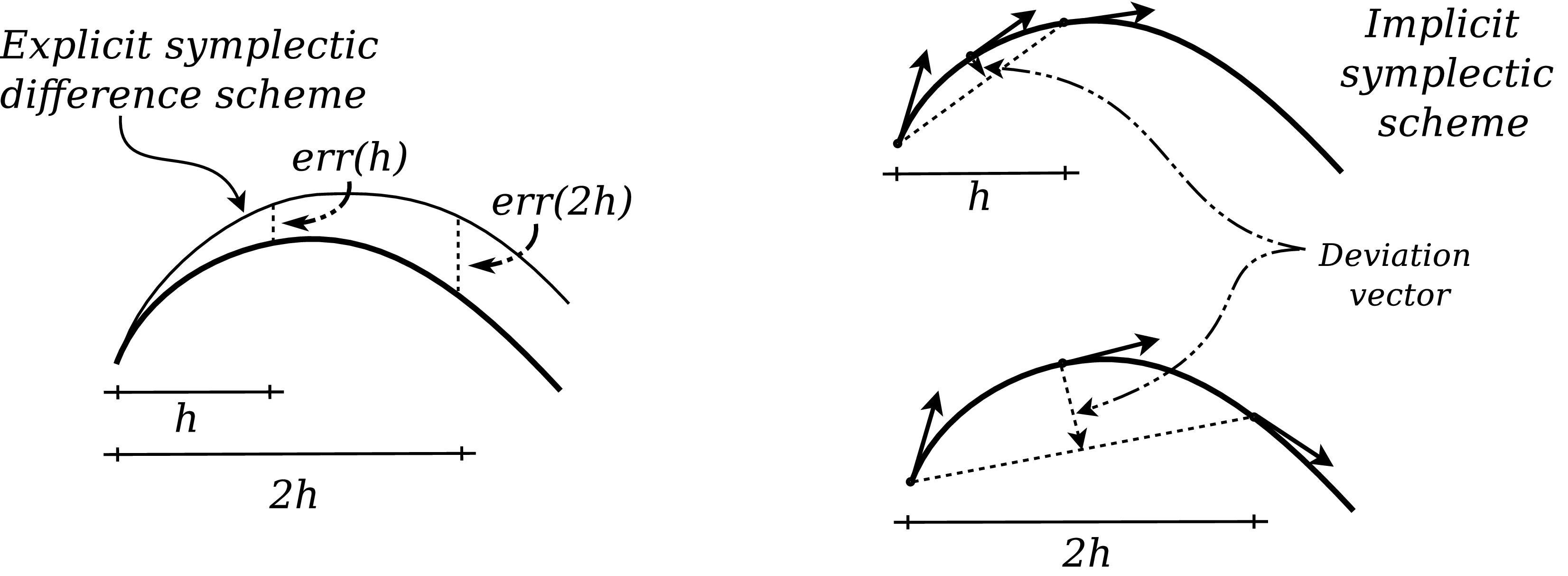}
 \caption{Explicit vs. implicit symplectic difference schemes. The deviation vector 
 in the implicit scheme introduced in this paper depends intrinsically on the time step $h$ and the
 Hamiltonian vector field $X_H$.}
 \label{fig:0}
\end{figure}

In order to find the right expression which gives a low order symplectic 
integrator as exact as possible using the method of Liouvillian forms, 
we construct a Hamiltonian isotopy between the continuous and the discrete
flows and we use the infinitesimal deformation of the isotopy for our symplectic
integrator. The vector field generated by this isotopy, is given in terms
of the original Hamiltonian vector field and is transversal to it.
By some classical relationships between Liouville and Hamiltonian vector fields,
with the Liouvillian forms we obtain the desired argument for our integrator. 

As in the precedent papers \cite{Jim15a,Jim15b}, we assume the reader 
is familiar with the terminology of differential geometry 
and vector bundles. For an introduction the reader is referred to \cite{AM78,LM87,MR91}.

\section{Hamiltonian and Liouville vector fields}

Consider a generic $2n$-dimensional manifold $M$ endowed with 
a symplectic form $\omega$, i.e. a non-degenerated, skew-symmetric, 
closed 2-form on $M$. 
The pair  $(M,\omega)$ is a 
\emph{symplectic manifold}. We say that it is \emph{exact} 
if the symplectic structure is exact, i.e., if there exists a primitive
1-form $\theta$ such that $\omega=d\theta$.
A \emph{Hamiltonian vector field} 
$X=X_H$ on $M$ is a vector field which satisfies 
$i_{X_H}\omega = -dH$
for a differentiable function $H:M\to\mathbb R$. 
The flow of a Hamiltonian vector field preserves the symplectic
form on $M$ which is characterized by the condition $\mL_{X_H}\omega=0$,
where $\mL_{X_H}\omega$ is the Lie derivative of $\omega$ along the integral curves 
of $X_H$. 
A \emph{Liouville vector field} $Z$ on a symplectic manifold 
$(M,\omega)$ is a vector field satisfying 
$\mathcal L_Z\omega=\omega$. Since $\omega$ is closed, the Lie 
derivative reduces to $\mL_Z\omega=d(i_Z\omega)=\omega$.
We write the 1-form as $\theta=i_Z\omega$, and we call it a \emph{Liouvillian 
form}\footnote{We use the term \emph{Liouville form} for the tautological 1-form
on the cotangent bundle $\pi:T^*Q\to Q$ given by $\pi^* \lambda\in T^*(T^*Q)$
and \emph{Liouvillian form} for the generic case $d\theta=\omega$.}
 \cite{LM87,Lib00}.
Several results and identities follow, in particular we have: 1) $i_Z\theta=0$, 2) $\theta=i_Zd\theta$ and
3) $\mL_Z\theta = \theta$, among many others.

The invariance of the symplectic form under the flow of Hamiltonian vector fields 
and the linearity of the Lie derivative show that Liouville vector fields
are invariant under the addition of Hamiltonian fields. Indeed, let $F:M\to\mathbb R$ be
a differentiable function, then  $Z+X_F$ is a Liouville 
vector field and by symplectic duality, $(\theta-dF)$ is a Liouvillian form. 

A symplectomorphism $\phi$ on an exact manifold $(M,d\theta)$ is 
called \emph{exact} with respect to the Liouvillian form $\theta$
if $\phi^*\theta-\theta=dF$ for a function $F:M\to\mathbb R$.

A \emph{symplectic isotopy} is a map $\phi_t:[0,1]\times M\to M: (t,q)\mapsto \phi_t(q)$
such that $\phi_t$ is a symplectic map for every $t$ and $\phi_0 = id$ and such that 
the vector field $X_t$ given by 
\begin{eqnarray}
   \frac{d}{dt}\phi_t = X_t\circ \phi_t,\qquad t\in[0,1], 
\end{eqnarray}
is symplectic. A \emph{Hamiltonian isotopy} is an exact symplectic isotopy,
it means that the vector field $X_t$ is Hamiltonian for every $t\in[0,1]$, i.e.
$i_{X_t}\omega=-dH_t$ for a time-dependent Hamiltonian function $H_t:[0,1]\times M\to \mathbb R$.

Some standard results in symplectic geometry relate the behaviour of a 
Liouvillian form under the flow of a symplectic and Hamiltonian isotopy. In particular
the following result is proved in \cite{MS17} for the case when $\theta$ is the tautological form.
\begin{proposition}
Let $(M,d\theta)$ be an exact symplectic manifold. An isotopy $\phi_t:[0,1]\times M\to M$ is symplectic if and only if 
$\alpha=\phi_t^*\theta -\theta$  is closed for every $t\in[0,1]$, and it is 
Hamiltonian if $dF_t=\phi_t^*\theta -\theta$ for a 1-parameter family of functions
$F_t:[0,1]\times M\to M$. 
\end{proposition}

\begin{remark}
 Note that if $\theta=p_0dq_0$ and $\phi:M\to M$ is a generic Hamiltonian symplectomorphism, 
 the pullback form $\phi^*\theta$ is not necesarily $\phi^*\theta=p_1dq_1$. 
 This fact, that is an usual trick 
 (see for example Remark 9.3.4 in \cite{MS17}), holds for very particular diffeomorphisms $\phi$
 called contact transformations or contactomorphisms defined on odd-dimensional 
 manifolds: either, codimension 1 submanifolds $\Sigma\subset M$, or the 
 $(2n+1)$-dimensional product
 $(M\times \mathbb R)$. 
 In fact, the method of Liouvillian forms is based on the fact that the 
 tautological form is not preserved under the flow of a generic Hamiltonian 
 flow. The idea is to find the 
 Liouvillian form whose variation depends on $H$.
\end{remark}

\section{Symplectic maps from Liouvillian forms}
We consider the results exposed in the
previous section for the construction of symplectic maps.
For this, we need to construct the geometrical framework
which is a classical procedure.

Define the product manifold of two 
copies of $(M,\omega)$ at times $t=0$ and $t=h$, which we denote by
$(M_1,\omega_1)$ and $(M_2,\omega_2)$. Assume that $(M_i,\omega_i)$
$i=1,2$, are diffeomorphic to cotangent bundles $M_i\cong T^*\mathcal Q_i$
where $\mathcal Q_i$, $i=1,2$, are configuration spaces of mechanical systems.
The canonical projections 
$\pi_i:\bP\to M_i$ for $i=1,2$ let us define a two-form $\omega_{\ominus}$ on $\bP$ by 
\begin{eqnarray}
    \omega_{\ominus} &=& \pi_1^*\omega_1 - \pi_2^*\omega_2.
    \label{eqn:def:sym}
\end{eqnarray}
The manifold $(\bP,\omega_{\ominus})$ becomes a symplectic manifold of dimension $4n$ 
\cite{LM87}. 

For any Liouvillian form $\theta$ on $\bP$, there exists a diffeomorphism 
$\Psi:\bP\to T^*(\mathcal Q_1\times \mathcal Q_2)$ such that 
$\theta = \Psi^*\theta_{\mathcal Q_1\times \mathcal Q_2}$. This diffeomorphism 
is symplectic and 
$(\bP,\mathcal Q_1\times\mathcal Q_2, \theta,\pi,\Psi)$ is a special
symplectic manifold on $\mathcal Q_1\times \mathcal Q_2$, where $\pi=\Psi^*\pi_{\mathcal Q_1\times \mathcal Q_2}$ 
\cite{Tul76,Tul77}.

Consider a function $F:\mathcal Q_1\times\mathcal Q_2\to \mathbb R$. The \emph{Lagrangian 
submanifold generated} by $F$ in the manifold $(\bP,d\theta)$ is defined by the equation
$\langle v, \theta\rangle = \langle T\pi(v), dF \rangle$ where $v\in T_p\bP$ and $\tau_{\bP}(v)=p$,
in the following way 
\begin{eqnarray}
    \Lambda = \left\{ p\in{\bP} | \pi(p)\in \mathcal Q_1\times \mathcal Q_2,\langle v, \theta\rangle =\langle T\pi(v),dF\rangle \right\}.
    \label{eqn:lag}
\end{eqnarray}

The submanifold $\Lambda$ is well defined since 
$\pi_{\mathcal Q_1\times\mathcal Q_2}:T^*(\mathcal Q_1\times\mathcal Q_2)\to\mathcal Q_1\times\mathcal Q_2$ is 
a submersion
\begin{eqnarray*}
    \xymatrix{
       v\in T \bP\ \ar[dr]_{\tau_{ \bP}} \ar[d]_{T\pi} & &  T^*  \bP\ni \theta \ar[dl]^{\pi_{ \bP}}  \\
       T(\mathcal Q_1\times\mathcal Q_2) \ar[dr]_{\tau_{\mathcal Q_1\times \mathcal Q_2}}  & p\in \bP \ \ar[d]_{\pi }  \ar[r]^{ \Psi\qquad } & 
              T^*(\mathcal Q_1\times\mathcal Q_2)\ni dF \ar[dl]^{\pi_{\mathcal Q_1\times \mathcal Q_2}}  \\
        & \mathcal Q_1\times\mathcal Q_2 &  
    }
\end{eqnarray*}

The preimage $\Psi^{-1}(\Lambda)\subset \bP$ is a Lagrangian submanifold in $\bP$. It 
corresponds to the graph $\Gamma_\phi$ of a symplectic map $\phi:M\to M$ by 
\begin{eqnarray}
   \Gamma_\phi = \left\{ (x,\phi(x))\in \bP \right\},
\end{eqnarray}
and it can be described by pulling-back the 1-form $\alpha=\Psi^*(dF)$ which is closed in $\bP$
but not necesarily exact.
We impose the condition that $\alpha$ be in addition, exact $\alpha=dS$, which implies
at the time, some restrictions on $\Phi\in Sp(\bP,\omega_\ominus)$. This fact is usually 
ignored since it is used to use $\Phi\equiv id$. We have two Lagrangian submanifolds 
$\Lambda\subset (T^*(Q_1\times Q_2), \omega_\oplus)$ and $\Gamma_\phi\subset(\bP,\omega_\ominus)$
defined by generating functions $dF$ and $dS=\Phi^*(dF)$, which concides with the 
restrictions of $\theta_{Q_1\times Q_2}|_\Lambda$ and $\theta|_{\Gamma_\phi}$ respectivelly.

The method of Liouvillian forms uses the (local) projection 
$\pi_N:U\subset\bP\to N$,
defined on a tubular neighborhood $U$ around $\Gamma_\phi$ by
\begin{eqnarray}
   \pi_N = J\circ(\pi_1 - \pi_2),
   \label{eqn:proj:def}
\end{eqnarray}
onto a $2n$-dimensional submanifold $N$.
This submanifold must behave like a symplectic submanifold
of $\bP$ and be related to the original manifold $(M,\omega)$.
This uses an additional symplectomorphism $\bP \leftrightarrow T^*M$ 
which corresponds to the well-known 1-to-1 correspondence between symplectic maps 
close to the identity with 1-forms close to the zero section in $T^*M$.
The projection (\ref{eqn:proj:def}) is in fact the projection
$\pi:U\subset\bP\to\Lambda$, it means that $\Lambda$ must be considered
as a 2n-dimensional submanifold in $\bP$ being symplectic for an 
alternative symplectic form $\tilde\omega$. 
Instead of constructing an additional 
special symplectic manifold, there is an easy way to deal with this 
extended framework.

We replace the geometry of the three 
symplectic manifolds $(T^*(Q_1\times Q_2),\omega_\oplus)$,
$(\bP,\omega_\ominus)$ and $(T^*M,\omega_{can})$ for a quaternionic
structure $\{I_{4n},\mI, \mJ, \mK\}$ on the product manifold $\bP=M_1\times M_2$
equiped with its natural Riemannian structure that we will denote by $\langle\cdot,\cdot\rangle$. 
It induces three different 
symplectic forms $\omega_\mI, \omega_\mJ,\omega_\mK$. Each symplectic form induces the geometry
of one of the 
previous symplectic manifolds. The $2n$-dimensional submanifold 
$\Lambda\in\bP$ which produces well defined symplectic maps for constructing symplectic
integrators must be Lagrangian for two of them and symplectic for the 
third. The details of this construction are given in \cite{Jim18}.

The projection given in (\ref{eqn:proj:def}) induces 
an intermediate point $\bar\bz=\rho(\bz_0,\bz_\tau)$, such that the 
implicit map given by
\begin{eqnarray}
   \bz_\tau = \bz_0 + \tau X_H\circ \rho(\bz_0,\bz_\tau).
   \label{eqn:map:liouvillian}
\end{eqnarray}
is symplectic if 
$\rho(\bz_0,\bz_\tau)$ satisfies the following two conditions
\begin{eqnarray}
   \frac{\partial \rho}{\partial \bz_0} + \frac{\partial \rho}{\partial \bz_\tau} = I_{2n},\quad {\rm and}\quad 
   \frac{\partial \rho}{\partial \bz_0} - \frac{\partial \rho}{\partial \bz_\tau} = b
   \label{eqn:cond}
\end{eqnarray}
where $b$ is a Hamiltonian matrix in $GL(2n,\mathbb R)$. We can write 
\begin{eqnarray}
   \rho(\bz_0,\bz_\tau) = \fiz(\bz_0 + \bz_\tau) + b(\bz_\tau -\bz_0),
\end{eqnarray}
moreover, we can substitute $b$ by $\tau b$ to have a symmetric integrator (see the details in \cite{Jim18,Jim15a}).
\begin{remark}
   This result was already obtained by Kang and his collegues, using the method
   of generating functions \cite{KQ10}. Their approach was mainly algebraic and only was 
   considered as a condition for obtaining an implicit symplectic map.
   We arrived to he same condition using Liouvillain forms and it gives a 
   geometrical interpretetion of the matrix $b$, as we will explain in the rest of this work.
\end{remark}

The matrix 
$b$ is related with the closed part of a Liouvillian form $\theta$ 
on $(M,\omega)$. Since all the computations are locally defined on open 
balls, by the Poincar\'e's lemma it corresponds to the exact part 
$df$ of the Liouvillian form, 
i.e. to the differential of a different generating function.
Moreover, using contact geometry, there is a way to associate a 
Liouvillian form to regular energy levels of a Hamiltonian function \cite{HZ12,MS17}.
Consequently, there is a well defined way to assign a $(1,1)$ 
tensor $\bb$ which generalizes the $b$ matrix for a prescribed 
regular energy level of a razonable Hamiltonian system $(M,\omega,X_H)$.

\section{Looking forward \emph{exact} symplectic integrators}
One way for minimizing the oscillations in a symplectic integrator is
measuring how much the discrete flow is far from the continuous
flow and correcting this deviation. We perform this task using
the results described in the previous sections applied to the flow of a Hamiltonian 
system $(M,\omega, X_H)$.
a geometrical construction which approximates the deviation
of the discretization for each $0\leq \tau<\tau_0$ for small $\tau_0$. 

Let fix the notation.
The flow of the Hamiltonian vector field $X_H$ will be denoted by
 $\varphi^t_H$, and it is solution of the 
$\dot \bz=X_H(\bz)$ with initial condition 
$\bz_0=\bz(0)=\varphi^0_H(\bz_0)$. We use alternatively the notation
$\bz(t)=\bz_t=\varphi_H^t(\bz_0)$. 

\begin{figure}[h!]
 \centering
 \includegraphics[scale=0.45]{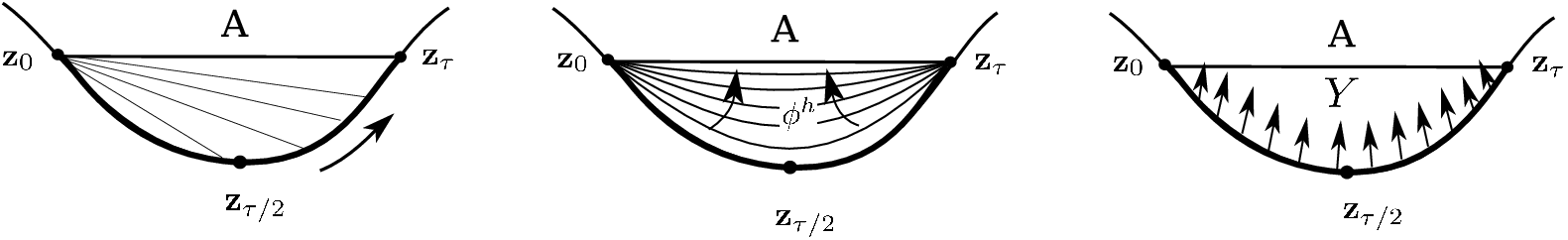}
 \caption{The isotopy $\psi_{h,\tau}^s$ and the vector field $Y$. Left: for every $0\leq \tau<\tau_0$
 we have a different segment ${\bf A}$ joining $\bz_0$ and $\bz_\tau$. Center:
 for every fixed $\tau$ the parameter $h$ gives the isotopy $\phi^h$ between the segment ${\bf A}$ and the Hamiltonian flow.
 Right: the local vector field $Y$ generated by the isotopy on the segment ${\bf A}$.}
 \label{fig:1}
\end{figure}

Let $0<\tau<\tau_0$ be a small value of $t$ and denote by  
${\bf A}$ the line segment joining $\bz_0= \bz(0)$ and
$\bz_\tau= \varphi_H^{\tau}(\bz_0)$.
The parameter $\tau$ represents the timestep of some discretization 
(left panel in Fig. \ref{fig:1}) and we consider that the mid-point is $t=\tau/2$
with value $\bz_{\tau/2}$. For small enough fixed values of $\tau_0$, we have a simple 
region enclosed by the segments ${\bf A}$ and $\bz([0,\tau])$,
that we can parameterize by two new real elements $s,h\in [1,0]$.
The parameter $h$ will determine an isotopy from $\varphi^t_H$ to ${\bf A}$
with fixed points $\bz_0$ and $\bz_\tau$. The parameter $s\in [0,1]$
will determine the curves joinning those fixed points (center panel 
in Fig. \ref{fig:1}). 
Since the segment ${\bf A}$ is given formally by the expression
${\bf A} = (1-s)\bz_0 + s\;\bz_\tau$, 
a first guess of this isotopy can be the following convex
parameterization
\begin{eqnarray}
  \psi_{h,\tau}^s(\bz_0)=(1-h)\varphi_H^{s\tau}(\bz_0) +h\left[(1-s)\bz_0 + s\;\bz_\tau\right], 
  \label{eqn:isotopy}
\end{eqnarray}
for every fixed $0\leq \tau<\tau_0$  and $h,s\in [0,1]$.

Unfortunatelly, just the segments given by the values $\tau=0$ and
$\tau=1$ in the parameterization (\ref{eqn:isotopy}) correspond to segments
of Hamiltonian flows.

We are looking for an isotopy $\phi^h:[0,1]\times M\to M$
that, written as a local map, must satisfies 
$$\phi^0 (\bz([0,\tau]))= \bz([0,\tau])\qquad {\rm and} \qquad 
\phi^1(\bz([0,\tau])) = {\bf A}$$ 
and such that the intermediate 
curves for $h\in (0,1)$ are also segments of Hamiltonian flows. We 
claim that this isotopy exists but the solution for our problem is, in fact,
much simpler. We will work with the Hamiltonian vector field $X_H$ and 
its pull-back by a local diffeomorphism $u\in {\rm Diff}(U)$ defined in an 
open set $U\subset M$ around the segment $\bz_{[0,\tau]}$. 

It is well-known that the pull-back of $X_H$ by a diffeomorphism 
$u:M\to M$ is given by 
$u^*X_H = (du)^{-1}\cdot X_H\circ u$ \cite{HZ12,Ste13}.
Also, we know that $u^*X_H$ is a Hamiltonian vector field 
if $u$ is a symplectomorphism. 
In the case of a one-parameter group of diffeomorphisms 
$u_t:\mathbb R\times M\to M$ the pull-back of $X_H$ is given by 
$$u_t^*X_H = (du_{-t})\cdot X_H\circ u_t.$$
Supose that the one parameter group $u_t$ has infinitesimal 
generator $Y$, then the pull-back of $X_H$ has infinitesimal
generator 
\begin{eqnarray} 
    X = \left.\left(\frac{d}{dt}  u_t^*X_H\right) \right |_{t=0} =  \mL_Y X_H
        = [X_H,Y].
    \label{eqn:gen}
\end{eqnarray}
Again, if $u_t$ is a one-parameter group of symplectomorphisms, then 
(\ref{eqn:gen}) is a Hamiltonian vector field with Hamiltonian function 
$\omega(X_H,Y) = -dH(Y)$. Moreover, if $Y=X_K$ is Hamiltonian with Hamiltonian
function $K$, then (\ref{eqn:gen}) has Hamiltonian function 
$\omega(X_H,X_K) = \{K,H\}$ which is the Poisson bracket of $H$ and $K$.
We have the classical relation $[X_H,X_K] = X_{\{K,H\}}$ which determines the 
Lie algebra isomorphism between functions and Hamiltonian vector fields.
All this applies on the whole manifold $(M,\omega)$ which we 
consider as the global case.

In the local case, there are local symmetries that cannot be 
extended to the global case. They are given by the flow of 
some Liouville vector fields $Z$ whose flow $\phi^h_Z$ satisfies 
$(\phi^h_Z)^*\omega = e^h\omega$. In the general case, the flow of 
$Z$ is non necesarily complet and we need to consider open neighborhoods $W$, 
big enough for including the source and target domains $W\subset \left(U\cap\phi^h_Z(U)\right)$, 
and small values of the parameter 
$h\in( -c,c)$ in the flow $\phi^h_Z$.

\begin{lemma}
   The pull-back of a Hamiltonian vector field $X_H$ under the 
   local flow $\phi^h_Z$ of a Liouville vector field $Z$ is a 
   (local) Hamiltonian vector field $X_{\mathcal A}$, with 
   Hamiltonian function $\mathcal A =  (\mL_Z H - H)$.
\end{lemma}
{\it Proof.} It is just the application of (\ref{eqn:gen}) for the case where 
the flow has a Liouville vector field as infinitesimal generator
\begin{eqnarray} 
    \left.\left(\frac{d}{dh}  \left(\phi^h_Z\right)^*X_H\right) \right |_{h=0} = \mL_Z X_H 
        = [X_H,Z].
    \label{eqn:par}
\end{eqnarray}
To prove that it is Hamiltonian, we check the contraction of
the vector field with the symplectic form $\omega$, 
indeed
$$i_{[X_H,Z]}\omega = \mL_{X_H}i_Z\omega - i_{Z}\mL_{X_H}\omega = \mL_{X_H}\theta 
  = -d(\mL_ZH - H).$$
Then we have a local function ${\mathcal A} = \mL_ZH - H$, and 
$[X_H,Z] = X_{\mathcal A}$ is a Hamiltonian vector field 
locally defined around the solution curves of $X_H$.
$\hfill\square$
\begin{remark}
   When the Liouville vector field is the vertical one
   $\sum_i p_i\frac{\partial}{\partial p_i}$, and the Hamiltonian 
   vector field is a natural mechanical system, the function $\mathcal A$
   corresponds to the Lagrangian function $L = 2T - H$, where $T$ is the 
   kinetic energy.
\end{remark}

\begin{remark}
   Previous lemma is a purely local result in contrast to the global 
   case considered by Theorem VI.2.8 in \cite{HLW10}. We conjecture that for Liouville fields
   that can be extended to the whole manifold, an additionnal first integral
   is concerned.
\end{remark}

The last discussion shows that we can use the flow of a Liouville vector field for constructing the 
local Hamiltonian isotropy connecting the segment $\bA$ with the Hamiltonian flow.   

\section{The geometrical meaning of the $\bb$ matrix}
As proved in \cite{KQ10} using generating functions, the map
\begin{eqnarray}
    \bz_\tau = \bz_0 + \tau X_H(\bar\bz), \qquad {\rm where} \qquad \bar\bz=\fiz(\bz_0+\bz_\tau) +\tau b(\bz_\tau-\bz_0)
    \label{eqn:map}
\end{eqnarray}
where $b\in \mathbb M_{2n\times 2n}(\mathbb R) $ is a Hamiltonian matrix, 
defines a symmetric, symplectic map for constructing a symplectic integrator.
In \cite{Jim15a,Jim18} this result was refined using Liouvillian forms,
where the matrix $b$ is generalized to a $(1,1)$ tensor $\bb$ on $M$, which corresponds to the 
closed (in fact to the exact) component of a Liouvillian form $\theta$ on $M$. 
Return to the Hamiltonian isotopy proposed in the last section and consider 
only the mid-point in the line segment $\bA$. The point $\bar\bz$
where the Hamiltonian vector field is evaluated in the implicit map 
(\ref{eqn:map}), is the image of the mid-point $\fiz(\bz_0+\bz_\tau)$
under the symplectic map $\phi=(I-2\bb_\tau)^{-1}(I+2\bb_\tau)$, where $\bb_\tau = \tau \bb$,
for enough small $\tau>0$. In this case, $\bb_\tau$
is non-exceptional and $\phi$ is well defined, moreover $\phi$ is close to the identity map. 

We construct the symplectic isotopy $\phi_h:[0,1]\to Sp(M,\omega)$ 
connecting the mid-point $\fiz(\bz_0+\bz_\tau)$ with $\bar\bz$
using the parameter $h$ by 
\begin{eqnarray}
   \phi_h = (I-h2\bb_\tau)^{-1}(I+h2\bb_\tau), \qquad h\in[0,1].
\end{eqnarray}
It satisfies $\phi_0=id$ and $\phi_1=\phi$, and it defines a symplectic map for each fixed $h\in[0,1]$.
Moreover, since it is close to the identity, it is a Hamiltonian isotopy for some 1-parameter family 
of Hamiltonian functions $H_t:[0,1]\times M\to \mathbb R$  \cite{MS17}.
Since $\phi$ is symplectic for every $\bb$ close to the zero tensor, the implicit map 
(\ref{eqn:map}) corresponds to the exact discretization of the flow of the 
Hamiltonian function $H\circ\phi:M\to \mathbb R$ known in the numerical community
as the ``sourrounding Hamiltonian''. To be more specific, the symplectic mid-point 
scheme exactly  integrates a ``surrounding Hamiltonian'' $\bar H = H\circ \varphi$
with equations of motion $\dot\zeta = X_{H\circ\varphi}(\zeta)$. 
Consequently, the map 
(\ref{eqn:map}) integrates exactly the system
\begin{eqnarray}
    \dot\zeta = X_{H\circ\varphi\circ \phi}(\zeta).
\end{eqnarray}

\begin{figure}[h!]
 \centering
 \includegraphics[scale=0.4]{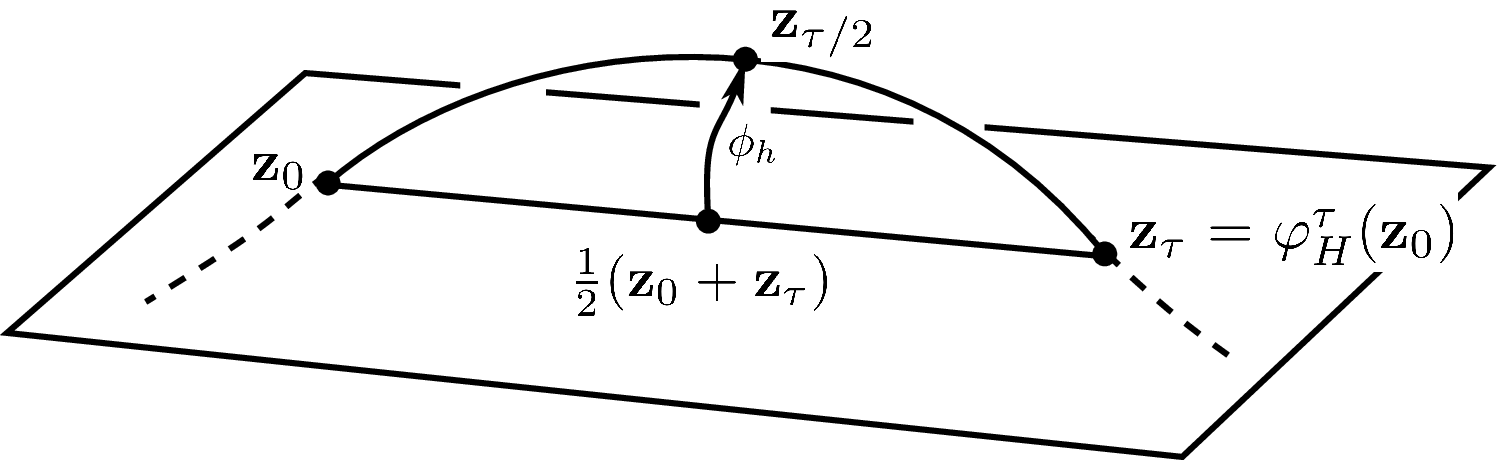}
 \caption{The Hamiltonian isotopy $\phi_h=(I-h2\bb_\tau)^{-1}(I+h2\bb_\tau)$, taking the mid-point to a point close to $\bz_{\tau/2}$.}
 \label{fig:2}
\end{figure}

The goal is to find the local symplectic map $\phi$ approximating $\varphi^{-1}$. Equivalently, we search for the $\bb$ tensor whose induced symplectic map $\phi$ takes the mid-point to some 
point on the solution we are integrating. If it is possible, $\phi^{-1}$ maps such a point on 
the mid-point, cancelling the numerical oscillations.
Before to propose some approaches for the search of the $\bb$ tensor, we will check the 
numerical algorithm.

\section{The symplectic integrator}

We consider the flow $\varphi^t_H$ of the Hamiltonian vector field 
$\dot \bz=X_H(\bz)$ and we integrate it in time from the initial condition 
$\bz_0$ 
\begin{eqnarray}
    \int_{0}^{\tau}\frac{d}{d t}\varphi_H^t(\bz_0)dt = \int_{0}^{\tau} X_H(\varphi_H^t(\bz_0))dt,\qquad \bz_0\in M, \quad
    \forall \tau\in I,
\end{eqnarray}
where $I = (-\tau_0,\tau_0)$ for small $\tau_0\in\mathbb R$.
Applying the fundamental theorem of calculus and
reparameterizing the time by $s(t)=\frac{1}{\tau}t$ we have  
\begin{eqnarray}
    \bz_\tau =\bz_0 + \tau \int_{0}^{1} X_H(\varphi_H^{\tau s}(\bz_0))ds,\qquad 
    \forall \tau\in I \subset\mathbb R,
    \label{eqn:par1}
\end{eqnarray}
which is the integral version of equation $\dot\bz=X_H(\bz)$). For small $\tau$, which is the case here,
expression (\ref{eqn:par1}) corresponds to the exponential map $\bz_\tau =e^{\tau X_H} \bz_0$.

The Cauchy-Lipschitz's theorem (a.k.a. Picard-Lindel\"{o}f's theorem) assures that a local solution 
for this equation always exists. Moreover, one way to approximate the value of $\bz_\tau$
is by means of Picard iterations. Given a first guess $\bz^1$ close to $\bz_\tau$ the 
iterative scheme 
\begin{eqnarray}
   \bz^{n+1} = \bz_0 + \int_0^\tau X_H(\bz^{n}(s))ds,
\end{eqnarray}
approximates the value $\bz_\tau$. The Picard-Lindel\"{o}f's theorem assures the 
convergence of this iterative process for small $\tau$. Note that for big 
values of $\tau$, the Lipschitz condition cannot be fulfilled.

Consider a symplectic integrator 
given by Picard iterations \cite{Jim15a}. Computing a first guess using an explicit 
scheme $\bz^1 = EulerA(\bz_0)$, we iterate
\begin{eqnarray}
    \bar\bz &=& \fiz(\bz_0 +\bz^n ) + \tau\bb(\bz^n-\bz_0)
    \label{eqn:point}\\
    \bz^{n+1} &=& \bz_0 + \tau X_H(\bar\bz).
\end{eqnarray}

This induces the following iterative algorithm:
\begin{center}
  \begin{tabular}{|rl|}
    \hline
    \multicolumn{2}{|l|}{
        {\bf Algorithm 1.} }\\
    \hline
        & $!$Setup the initial guess for $\bz_{n+1}$\\
     1: & $\bz^{[0]}=\bz_n + \tau X_H(\bz_n)$\\
     2: & for $j=0:\kappa$ do\\
        & $\quad$ $!$ compute the tensor $\bb$ at the mid-point \\
     3: & $\quad b =  \bb \left( \fiz( \bz_n + \bz^{[j]}) \right)$ \\
        & $\quad$ $!$ compute the point $\bar \bz$ \\
     4: & $\quad \bar{\bz} =  \fiz(\bz_n + \bz^{[j]})  + \tau b ( \bz^{[j]}-\bz_n )$ \\
        & $\quad$ $!$ refine the guess\\
     5: & $\quad$ $\bz^{[j+1]}= \bz_n + \tau X_H(\bar \bz)$ \\
     6: & end for\\
     7: & $\bz_{n+1}=\bz^{[\kappa]}$\\
    \hline
  \end{tabular}
\end{center}

The challenge is to find the way to compute, for a prescribed
Hamiltonian system, a good guess of $\bb$  depending
in addition on the parameter $\tau$.
We can consider the  value of $\bar\bz$ given in (\ref{eqn:point}) as a first order approximation in $\tau$,
and develop the tensor $\hat\bb$ as a series in even powers of $\tau$
as follows
\begin{eqnarray}
    \hat\bb =  \bb_0 + \tau^2 \bb_2 + \tau^4 \bb_4 + \cdots 
    \label{eqn:bhat}
\end{eqnarray}
where $\bb_j$, $j=2k$ are symmetric, 
Hamiltonian $(1,1)$-tensors. 
Inserting (\ref{eqn:bhat}) in (\ref{eqn:point}) the symmetry 
$(\tau,\bz_0,\bz_\tau) \mapsto (-\tau,\bz_\tau,\bz_0)$ is preserved
and the integrator preserves symmetry and symplecticity.

\subsection{Looking forward the associated Liouvillian form}

We have shown that the path that takes the mid-point to the point $\bar \bz$
given in (\ref{eqn:map})
is a Hamiltonian isotopy. This isotopy is attached to the 
Liouvillian form $\theta$ which defines the map. On the other hand, in contrast 
to the symplectic form $\omega$
which is preserved by the Hamiltonian flow $\varphi^t_H$, a Liouvillian form is not preserved, 
but it produces a Hamiltonian isotopy $dF_t:=(\varphi^t_H)^*\theta - \theta$ which is related 
to the previous isotopy but they are not the same. 
We want to find a Liouvillian form induced by the geometry
of the Hamiltonian system $(M,\omega, X_H)$, such that its Liouville vector field determines
an infinitesimal generator on a prescribed solution which sends it into another local solution 
but only in a tubular neighborhood.

The solution to this problem is adapted from an equivalent problem in the interface
of contact and symplectic topologies. In the terminology used by McDuff and Salamon \cite{MS17}
it concerns the \emph{internal symplectization} of a contact manifold. This procedure
in addition, imposes a constraint concerning the lenght of the segment of Hamiltonian flow 
where the method works. In fact, this constraint comes from Gromov's non-squeezing theorem
\cite{Gro86} and it is related to the symplectic width of the ``symplectized manifold''.
This explains why Ge and Marsden's lemma \cite{GM88} on the reparameterization 
of the Hamiltonian flow is not a sufficient condition (among others assumptions) 
for claiming the non existence of energy preserving symplectic integrators. 

The procedure to find a Liouvillian form for $(M,\omega,X_H)$ works for regular solutions,
i.e. for solutions belonging to a regular level hypersurface. The interested readers
are refered to \cite{HZ12,MS17,Wei71} for the generic construction, and \cite{Jim18}
for the procedure adapted to a prescribed Hamiltonian system. In this paper we will 
only sketch the global procedure. It can be explained in two big steps. 

The first step consists in to define 
a contact structure on the regular hypersurface fixed by the initial condition $\bz_0$,
leading to a contact manifold embedded in $(M,\omega)$. For this, fix the level hypersurface
using the initial condition $\bz_0\in M$ with regular value $H(\bz_0)=h_0\in \mathbb R$. 
The set $\Sigma_{h}=H^{-1}(h_0)$ is a smooth submanifold of codimension 1 
by Saard's theorem. Select a 
Liouville vector field $Z$ on $(M,\omega)$ which is transversal to $\Sigma_h$
and regular in a tubular neighborhood $U\subset M$ around $\Sigma_h$. 
Consider $\Sigma_h\subset U$ as an embedding $\jmath:\Sigma_h\to U$, and define
the linear form $\alpha=\jmath^*\theta$ on $\Sigma_h$, which is the pullback of the 
Liouvillian form $\theta = i_Z\omega$ to $\Sigma_h$. 
We consider the distribution 
$\mathcal D=\ker(\alpha)$ which endows $\Sigma_h$ with a contact 
structure. $(\Sigma_h,\mathcal D)$ becomes a contact manifold 
with contact form $\alpha$. Finally, we define a Reeb field for $\alpha$ from the Hamiltonian
vector field restricted to $\Sigma_h$ given by 
$R_\alpha = \frac{1}{\alpha(X_H|_\Sigma)}X_H|_\Sigma$.
Note that $R_\alpha$ is a rescaling of the Hamiltonian vector field and it depends on the 
selected vector field $Z$.

The second step is the internal symplectization which can be splitted in two parts:
1) the external symplectization, mapping 
$(\Sigma_h,\alpha)\to(\Sigma_h\times\mathbb R,d(e^s\alpha))$ and 2) the embedding of
a slice around $\Sigma_h$ into the tubular neighborhood $U$, in this way
$(\Sigma_h\times (-c,c), d(e^s\alpha))\to (U,d(\theta_{H,h}))$. Note that it is not 
necessary that $Z$, $\theta$ or $\theta_{H,h}$ be global forms, it suffies 
their regularity in the tubular neighborhood $\Sigma_h\subset U$.
The difficulty in passing from the external to the internal symplectization is the 
construction of the Liouvillian form $\theta_{H,h}$ from 
$e^s\alpha$. It is a classical procedure obtained by using Weinstein's proof of Darboux's theorem,
Moser's trick and the homotopy lemma
\cite{MS17,LM87,HZ12,Wei71}. 

The following expression is proved in \cite{Jim18}: The Liouvillian form 
$\theta_{H,h}$ associated to the Hamiltonian $H:M\to \mathbb R$ at the 
hypersurface level $h=H(\bz_0)$ is given by 
\begin{eqnarray}
   \theta_{H,h} = e^{h-H}\pi^*\alpha + \Phi(\varpi)
   \label{eqn:liouvillian}
\end{eqnarray}
where 
\begin{eqnarray*}
   \varpi = \left(1-e^{h-H}\right)\omega +  e^{h-H}dH\wedge \theta 
\end{eqnarray*}
is a closed 2-form which vanishes on $\Sigma_h$, and  $\Phi(\varpi)$ is its cochain homotopy
\begin{eqnarray}
    \Phi(\varpi)=\int_0^1(\psi^*_s(i_{X_s}\varpi))ds
    \qquad{\rm where}\qquad X_s=\frac{s}{\|\nabla H\|^2}\nabla H,
\end{eqnarray}
and $\psi_s$ is the flow of the rescaled gradient $X_s$.

Once the Liouvillian form $\theta_{H,h}$ was computed, we extract the symmetric 
part which belongs to the kernel of the differential $d\theta_{H,h}=\omega$. Just for simplicity, we 
consider local coordinates $\{z_i\}_{i=1}^{2n}$ on $(M,\omega)$,  related to Darboux's coordinates by
$(q_i,p_i)=(z_i,z_{n+i}), i=1,...,n$. The Liouvillian 
form has a local expression in these coordinates by 
$\theta_{H,h}=\sum_i \alpha_i(z)dz_i$, where $\alpha_i:M\to\mathbb R$ are smooth 
functions. The closed part of $\theta_{H,h}$
is given by the symmetric matrix $S=(S_{ij})$ with expression
\begin{eqnarray}
S_{ij} = \fiz\left( \frac{\partial \alpha_i}{\partial z_j} + \frac{\partial \alpha_j}{\partial z_i}\right).
\end{eqnarray}

Finally we obtain a tensor $\tilde\bb$ which contains the information of the Hamiltonian 
flow at the energy level $H(\bz_0)=h$ by $\tilde\bb = JS$,\footnote{Other 
possibilities are $-JS$, $SJ$ and $-SJ$, since all of them are Hamiltonian.}
where $J$ is the 
complex structure associated to $\omega$. Inserting $\tilde\bb$ into 
(\ref{eqn:point}) we obtain a symplectic integrator adapted for simulating the 
flow $\varphi^t_H$ of the Hamiltonian vector field $X_H$ with initial condition 
$\bz_0$.

The tensor $\bb$ which minimizes the oscillations is related to 
$\theta_{H,h}=i_{(Z_H,h)}\omega$,
since the Liouville vector field $Z_{H,h}$ is the 
infinitesimal generator of the Hamiltonian isotopy connecting the continuous flow
with the mid-point numerical approximation. This relation between $\bb$ and 
$\theta_{H,h}$ can be highly non-linear. In a future it can be interesting 
to study this problem from the variational point of view.

\section{Conclusions and perspectives}
In this paper we refined the numerical scheme introduced in \cite{Jim15a} and 
we collected a series of results to give a full geometrical explanation of the 
method of Liouvillian forms with application to symplectic integration. 
The geometric approach gives an intuitive framework for understanding
the oscillatory behaviour of the numerical solution produced by a symplectic 
integrator when simulating Hamiltonian dynamics. A symplectic integrator defines 
intrinsically a Liouville vector field $Z$ and visceversa. The oscillations 
correspond to the projection of $Z$ on the gradient vector field $\nabla H$.
This is $$\langle Z,\nabla H\rangle = dH(Z)=\omega(Z,X_H)=\theta(X_H)$$
At first sight the method 
looks cryptic and abstract since the technique for finding the Liouvillian form 
associated to a prescribed Hamiltonian is difficult to visualize. However, this method
shows that there is no local obstruction for approximating the separation
of the numerical solution with respect to the continuous solution.  

We obtain a framework which extends the method of generating functions 
giving an algorithmic way for constructing a symplectic map for approximating
the flow of (almost) any generic, natural and classical Hamiltonian system.
The use of a quaternionic structure on the product symplectic manifold
simplifies the framework of special symplectic manifolds and gives a 
geometrical explanation to the Hamiltonian matrix $b$, first studied by Feng Kang 
as a condition for constructing implicit symplectic maps \cite{KQ10}.
The quaternionic structure shows the relation between four objects: 
1) the $b$ matrix which extends to the $\bb$ tensor in this framework,
2) the Liouvillian form where the symmetric part $S$ of its diferential induces 
$\bb=JS$, 3) the element 
$\bar\bz=\fiz(\bz_0+\bz_\tau)+\bb(\bz_\tau,\bz_0)$ which is interpreted 
as a tangent vector to the Lagrangian submanifold containing the flow, 
and 4) the symplectic map $\phi=(I-2\bb)^{-1}(1+2\bb)$ which 
is the (symplectic) Cayley transformation of $\bb$  \cite{Jim19a}. Moreover,
if the ``surrounding Hamiltonian'' of the mid-point rule is $\bar H = H\circ \varphi$ then the ``surrounding Hamiltonian'' of this symplectic integrator 
is $\hat H=H\circ\varphi\circ\phi$. The challenge is to approximate 
the map $\varphi^{-1}$ in order to 
approximate the original function $H$.

Once the relationship between the Liouvillian form, the symplectic map
and the symplectic integrator is given, we search for a suitable tensor 
$\bb$ for inserting in the numerical scheme. Fortunatelly, there is a rich 
theory for the search of closed characteristics on compact contact type manifolds
concerning a celebrated conjecture stated by Weinstein \cite{HZ12}. One of the 
main tools for solving this conjecture is the construction of a Liouville vector 
field, which is transversal to the contact type manifold at every point. We 
adapt this technique for a prescribed Hamiltonian system and we construct the 
Liouvillian form $\theta_{H,h}$ in a tubular neighborhood around the level 
hypersurface which 
contains the initial condition. This procedure is developed with all the 
details in \cite{Jim18}. 
This closes the loop relating the Hamiltonian system (with a given 
initial condition $\bz_0$), the Liouvillian form, the $\bb$ tensor and the 
symplectic integrator.
\begin{remark}
Note that the Liouvillian form 
$\theta_{H,h}$ depends on the Hamiltonian $H$ and, must important, 
it depends explicitly on the value $h=H(\bz_0)$. For Hamiltonians with no 
other first integral and for chaotic systems, different values $h$ and $h+\epsilon$,
for small $0<\epsilon\ll 1$, produce different Liouvillian forms.
\end{remark}

A systematic study on numerical techniques for approximating the tensor 
$\bb$, must be put in practice. In addition, it is necessary to search for 
a practical way of computing $\bb$ without computing the integral expression
of the cochain homotopy in the Liouvillian form (\ref{eqn:liouvillian}). 
An alternative that we will study in the future is the approximation of the 
tensor $\bb$ using variational methods.

\section*{Acknowledgements}
The author thanks J.P. Vilotte and B. Romanowicz for their support and constructive
criticism on this work. 
Special thanks to Profr. Robert McLachlan for signaling to me a problem with a 
previous interpretation of the object modifying the mid-point in the algorithm 
of the first version of this paper.
This research was developed with support from the \emph{Fondation du Coll\`ege de 
France} and \emph{Total} under the research convention PU14150472, as well as the ERC Advanced Grant 
WAVETOMO, RCN 99285, Subpanel PE10 in the F7 framework.

\end{document}